\documentclass[12pt]{article}
\usepackage{amssymb}
\begin{document}
$\;$\\[-27mm]
\begin{center}
{\Large\bf Non-Lie subgroups in Lie groups over local\\[2mm] fields of positive characteristic}\\[7mm]
{\bf Helge Gl\"{o}ckner\footnote{Supported by Deutsche Forschungsgemeinschaft,
project GL 357/10-1.}}\vspace{4mm}
\end{center}
\begin{abstract}\vspace{.3mm}\noindent
By Cartan's Theorem, every closed subgroup $H$ of a real (or $p$-adic)
Lie group $G$ is a Lie subgroup.
For Lie groups over
a local field ${{\mathbb K}}$ of positive characteristic,
the analogous conclusion is known to be wrong.
We show more:
There exists a ${{\mathbb K}}$-analytic Lie group $G$ and a non-discrete, compact subgroup $H$
such that, for every ${{\mathbb K}}$-analytic manifold~$M$,
every ${{\mathbb K}}$-analytic map $f\colon M\to G$ with $f(M)\subseteq H$
is locally constant. In particular,
the set~$H$ does not admit
a non-discrete ${{\mathbb K}}$-analytic manifold structure
which makes the inclusion
of $H$ into $G$ a ${{\mathbb K}}$-analytic map.
We can achieve that, moreover,
$H$ does not admit a ${{\mathbb K}}$-analytic Lie group
structure compatible with the topological group structure
induced by~$G$ on~$H$.\vspace{3mm}
\end{abstract}
{\footnotesize {\em Classification}:
22E20 (primary); 
22E35, 
22E50, 
32P05 (secondary). \\[3mm]
%
%
%
%
}
\section{Introduction and statement of results}
If ${{\mathbb K}}$ is the field ${{\mathbb Q}}_p$ of $p$-adic numbers or the field ${{\mathbb R}}$
of real numbers,
then Cartan's Theorem
ensures that every closed subgroup~$H$
of a ${{\mathbb K}}$-analytic Lie group~$G$
is a submanifold
and hence a Lie subgroup
(see the corollary to Theorem~1 in \cite[Part~II, Chapter~V, \S9]{Ser});
viz., for each $x\in H$, there exists
a ${{\mathbb K}}$-analytic diffeomorphism $\phi$
from an open $x$-neighbourhood $U\subseteq G$
onto an open subset $V$ of $E\cong{{\mathbb K}}^n$
such that
\[
\phi(U\cap H)=V\cap F
\]
for some vector subspace $F\subseteq E$.
It is well known that the conclusion
becomes invalid for Lie groups
over the field ${{\mathbb C}}$ of complex numbers
or a local field ${{\mathbb K}}$
of positive characteristic,
i.e., a field ${{\mathbb K}}={{{\mathbb F}}\mkern+.2mu(\mkern-2.3mu(X)\mkern-2.3mu)}$
of formal Laurent series over a finite field~${{\mathbb F}}$
(see \cite{Wei} for basic information concerning such fields).
In fact, ${{\mathbb R}}$ is a closed subgroup
of the $1$-dimensional complex Lie group $({{\mathbb C}},+)$
which cannot be a Lie subgroup
as it is neither open nor discrete.
For the same reason, the closed subgroup
\[
{{\textstyle{{\mathbb F}}\mkern-1.7mu \left[\mkern-1mu \left[X^2\mkern+1.3mu\right]\mkern-1mu \right]}}
:=\left\{\sum_{k=0}^\infty a_kX^{2k}\colon (a_k)_{k\in{{\mathbb N}}_0}\in{{\mathbb F}}^{{{\mathbb N}}_0}\right\}
\]
is not a Lie subgroup of the $1$-dimensional
${{\mathbb K}}$-analytic Lie group $({{\mathbb K}},+)$, if ${{\mathbb K}}={{{\mathbb F}}\mkern+.2mu(\mkern-2.3mu(X)\mkern-2.3mu)}$.
Our goal is to describe subgroups of Lie groups
over ${{\mathbb K}}={{{\mathbb F}}\mkern+.2mu(\mkern-2.3mu(X)\mkern-2.3mu)}$
with more special properties.\\[2.3mm]
Recall that a Lie group over a locally compact field
of characteristic~$0$ induces a Lie group structure
on each subgroup (see Definition~3
in \cite[Chapter~III, \S4.5]{Bou}).
In the example of ${{\mathbb R}}\subseteq{{\mathbb C}}$ just recalled,
the induced structure is discrete,
as observed in \cite[Chapter~III, \S4.5]{Bou}.\\[2.3mm]
We show that the non-discrete, compact subgroup
\begin{equation}\label{theone}
H:=\left\{\sum_{k=0}^\infty a_kX^{2^k}\colon (a_k)_{k\in{{\mathbb N}}_0}\in{{\mathbb F}}^{{{\mathbb N}}_0}\right\}
\end{equation}
of $({{\mathbb F}}(\!(X)\!),+)$
has analogous properties.
With terminology as in \S\ref{secprel},
we have:\\[2.5mm]
{\bf Theorem 1.1.}
\emph{Let ${{\mathbb K}}:={{{\mathbb F}}\mkern+.2mu(\mkern-2.3mu(X)\mkern-2.3mu)}$ for a finite field~${{\mathbb F}}$.
The subgroup $H$ of the $1$-dimensional
${{\mathbb K}}$-analytic Lie group $G:=({{\mathbb K}},+)$
described in {\rm(\ref{theone})}
has the following properties}:
\begin{itemize}
\item[\rm(a)]
\emph{$H$ is an
initial Lie subgroup of $G$.
The initial Lie subgroup topology on~$H$ is the discrete
topology.}
\item[\rm(b)]
\emph{For each ${{\mathbb K}}$-analytic manifold structure on the set~$H$
making the inclusion map $H\to G$ a ${{\mathbb K}}$-analytic map,
the topology on the manifold~$H$ must be discrete.}\vspace{.5mm}
\end{itemize}
\noindent
In particular, $H$ is not a submanifold of~$G$.\\[2.3mm]
In Section~\ref{secprel}, we introduce certain
\emph{constancy-enforcing} subsets~$N$ of
a ${{\mathbb K}}$-analytic manifold~$M$.
These can be regarded as initial submanifolds
such that~$N$ as a manifold is endowed with the discrete
topology (see Remark~2.5).
Using this terminology,
Theorem~1.1\,(a) requires that~$H$
be a constancy-enforcing subset of~$G$.
Part~(b) follows from~(a).\\[2.3mm]
Note that $H$ from~(\ref{theone})
can be given a $1$-dimensional ${{\mathbb K}}$-analytic Lie group structure
making the bijection
\[
{{\textstyle{{\mathbb F}}\mkern-1.7mu \left[\mkern-1mu \left[X\mkern+1mu\right]\mkern-1mu \right]}}\to H,\quad \sum_{k=0}^\infty a_kX^k\to \sum_{k=0}^\infty a_kX^{2^k}
\]
an isomorphism of ${{\mathbb K}}$-analytic Lie groups,
where ${{\textstyle{{\mathbb F}}\mkern-1.7mu \left[\mkern-1mu \left[X\mkern+1mu\right]\mkern-1mu \right]}}\subseteq{{{\mathbb F}}\mkern+.2mu(\mkern-2.3mu(X)\mkern-2.3mu)}$
is the compact open subring of formal power series.
This ${{\mathbb K}}$-analytic Lie group structure on~$H$
is compatible with the given compact topology on~$H$,
induced by $G$;
it does not make the
inclusion map $H\to G$ a ${{\mathbb K}}$-analytic mapping,
by~(b).\\[2.3mm]
In our second example,
a manifold structure compatible with the
induced topology does not exist.
For ${{\mathbb K}}={{{\mathbb F}}\mkern+.2mu(\mkern-2.3mu(X)\mkern-2.3mu)}$, we consider
the unit group
$({{\mathbb K}}^\times,\cdot)$ of non-zero field elements,
which is a $1$-dimensional ${{\mathbb K}}$-analytic Lie group.
It is well known that ${{\mathbb K}}^\times$ contains a compact
open subgroup isomorphic to a countable
power $({{\mathbb Z}}_p)^{{\mathbb N}}$ of the compact group ${{\mathbb Z}}_p$
of $p$-adic integers,
where $p$ is the characteristic of~${{\mathbb F}}$
(cf.\
\cite[Chapter~II, \S3, Proposition~10]{Wei}). Hence ${{\mathbb K}}^\times$ has a compact subgroup
$H$ which is isomorphic to ${{\mathbb Z}}_p$
as a topological group.
We can take $H$ as before and $G:={{\mathbb K}}^\times$
in the following result.\\[2.5mm]
{\bf Theorem~1.2.}
\emph{Let $G$ be a Lie group over
${{\mathbb K}} \! =\! {{{\mathbb F}}\mkern+.2mu(\mkern-2.3mu(X)\mkern-2.3mu)}$
and $H\!\subseteq \!G$ be a compact subgroup
such that $H\cong{{\mathbb Z}}_p$
as a topological group.
Then $H$ does not admit a ${{\mathbb K}}$-analytic
Lie group structure
compatible with the topology induced by~$G$.}\\[2.5mm]
The pathologies described
in Theorem~1.1
and Theorem~1.2
can occur simultaneously,
as the following example shows.\\[2.5mm]
{\bf Example 1.3.}
Let $p$ be a prime,
${{\mathbb F}}$ be the field with $p$ elements
and $2\leq \ell\in {{\mathbb N}}$
such that $p$ and $\ell$ are coprime.
Then the closed subgroup~$H$
generated by $1+X^\ell$ in the multiplicative group
$G:={{{\mathbb F}}\mkern+.2mu(\mkern-2.3mu(X)\mkern-2.3mu)}^\times$
is isomorphic to ${{\mathbb Z}}_p$,
whence it does not admit a Lie group structure
over ${{\mathbb K}}:={{{\mathbb F}}\mkern+.2mu(\mkern-2.3mu(X)\mkern-2.3mu)}$
compatible with its topology
(by Theorem~1.2).
Note that $H\subseteq {{\textstyle{{\mathbb F}}\mkern-2.2mu \left[\mkern-2.8mu \left[X^\ell\mkern+1.5mu\right]\mkern-2.8mu \right]}}$.
Thus~$H$ is a constancy-enforcing subset
of~$G$
(and hence an initial Lie subgroup
when endowed with the discrete
topology), by
the following theorem.\\[2.5mm]
{\bf Theorem~1.4.}
\emph{Let $p$ be a prime,
${{\mathbb F}}$ be the field of order~$p$
and $2\leq \ell\in {{\mathbb N}}$
with $p$ and $\ell$ coprime.
Then ${{\textstyle{{\mathbb F}}\mkern-2.2mu \left[\mkern-2.8mu \left[X^\ell\mkern+1.5mu\right]\mkern-2.8mu \right]}}$
is a constancy-enforcing subset of
${{{\mathbb F}}\mkern+.2mu(\mkern-2.3mu(X)\mkern-2.3mu)}$.}\\[2.5mm]
{\bf Remark~1.5.}
If $H$ is as in Example~1.3,
then also every subgroup $S\subseteq H$ is constancy-enforcing.
Moreover, $S$ does not admit a non-discrete ${{\mathbb K}}$-analytic Lie
group structure making the inclusion map $S\to G$ continuous.
[Other\-wise, $S$ would have a non-trivial compact open subgroup~$K$.
Then $K$ would be a non-trivial closed subgroup of $H\cong {{\mathbb Z}}_p$
and thus $K\cong{{\mathbb Z}}_p$,
contrary to Theorem~1.2].\\[1.2mm]
\emph{General conventions.}
We write ${{\mathbb N}}:=\{1,2,\ldots\}$ and ${{\mathbb N}}_0
:={{\mathbb N}}\cup\{0\}$.\\[2.4mm]
\noindent
{\bf Acknowledgement.}
The author thanks C.R.E. Raja
(Indian Statistical Institute, Bangalore)
for questions and discussions
which inspired the work.
The referee's comments helped to improve
the presentation.
\section{\hspace*{-3.8mm}Initial\hspace*{-.2mm} submanifolds\hspace*{-.2mm}
and\hspace*{-.2mm} constancy-enforcing\\
\hspace*{-3.6mm}subsets}\label{secprel}
Recall that a topological field~${{\mathbb K}}$ is called
a \emph{local field} if it is non-discrete
and locally compact (see \cite{Wei}).
We fix an absolute value $|\cdot|$ on~${{\mathbb K}}$
which defines its topology.
For the basic theory of ${{\mathbb K}}$-analytic mappings
between open subsets of finite-dimensional
${{\mathbb K}}$-vector spaces and the corresponding
${{\mathbb K}}$-analytic manifolds and Lie groups modeled on ${{\mathbb K}}^m$,
we refer to \cite{Ser}
(cf.\ also \cite{FAS} and \cite{Bou}).
The following well-known concept is useful.\\[2.5mm]
{\bf Definition 2.1.}
Let $M$ be a ${{\mathbb K}}$-analytic manifold over a local field~${{\mathbb K}}$.
A subset $N\subseteq M$
is called an \emph{initial submanifold}
if there exists a ${{\mathbb K}}$-analytic manifold
structure on~$N$ with the following property:
For each ${{\mathbb K}}$-analytic manifold~$Z$,
each map $f\colon Z\to N$
is analytic if and only if $\iota_N\circ f\colon Z\to M$ is analytic,
where $\iota_N\colon N\to M$, $x\mapsto x$ is the inclusion map.\\[2.5mm]
{\bf Remark~2.2.}
\begin{itemize}
\item[(a)]
If $N\subseteq M$ is an initial submanifold,
then $\iota_N\colon N\to M$ is analytic.
In fact, $\mbox{id}_N\colon N\to N$
is analytic and hence also $\iota_N\circ \mbox{id}_N=\iota_N$.
\item[(b)]
The ${{\mathbb K}}$-analytic manifold structure
on an initial submanifold is uniquely determined.\\[2mm]
In fact, let us write $N_1$ and $N_2$ for~$N$, endowed
with ${{\mathbb K}}$-analytic manifold
structures as described in Definition~2.1.
Write $\iota_{N_j}$ for the inclusion map $N_j\to M$, for $j\in\{1,2\}$.
Let $f\colon N_1\to N_2$ be the identity map, $x\mapsto x$.
Since $\iota_{N_2}\circ f=\iota_{N_1}$ is analytic,
the map $f$ is analytic.
Likewise, the analyticity of $\iota_{N_1}\circ f^{-1}=\iota_{N_2}$
implies that $f^{-1}$ is analytic. Thus $f\colon N_1\to N_2$
is an analytic diffeomorphism.
\item[(c)]
Since $\iota_N$ is ${{\mathbb K}}$-analytic and hence
continuous, the topology on an initial submanifold
(as a manifold) is finer than the topology induced by~$M$.
It can be properly finer.\vspace{.5mm}
\end{itemize}
{\bf Definition 2.3.}
If $G$ is a ${{\mathbb K}}$-analytic Lie group,
we call a subgroup $H\subseteq G$
an \emph{initial Lie subgroup}
if it is an initial submanifold.\\[2.5mm]
If $m_H\colon H\times H\to H$ is the group multiplication of~$H$
and $m_G$ the group multiplication of~$G$,
then $\iota_H\circ m_H=m_G\circ (\iota_H\times\iota_H)$
is ${{\mathbb K}}$-analytic, entailing that $m_H$ is ${{\mathbb K}}$-analytic.
Likewise, the map $H\to H$, $x\mapsto x^{-1}$
is ${{\mathbb K}}$-analytic, and thus $H$ is a ${{\mathbb K}}$-analytic
Lie group in its initial manifold structure.\\[2.3mm]
We consider a related concept.\\[2.5mm]
{\bf Definition 2.4.}
Let $M$ be a ${{\mathbb K}}$-analytic manifold over a local field ${{\mathbb K}}$.
We say that a subset $N\subseteq M$ is
\emph{constancy-enforcing}
if,
for each ${{\mathbb K}}$-analytic manifold~$Z$,
each ${{\mathbb K}}$-analytic map $f\colon Z\to M$
with image $f(Z)\subseteq N$ is locally constant.\\[2.5mm]
{\bf Remark 2.5.}
Note that $N$ is constancy-enforcing if and only
if $N$ is an initial submanifold
of~$M$ and the topology underlying the
initial submanifold structure
is the discrete topology.
For our ends, the shorter phrase is valuable.\\[2.5mm]
Let us write ${{\mathbb D}}:=\{z\in{{\mathbb K}}\colon |z|<1\}$
for the open unit ball.\\[2.5mm]
{\bf Lemma~2.6.}
\emph{A subset $N\subseteq M$ is constancy-enforcing
if and only if each ${{\mathbb K}}$-analytic map
$f\colon{{\mathbb D}}\to M$ with $f({{\mathbb D}})\subseteq N$
is constant on some $0$-neighbourhood.}\\[2.5mm]
{\bf Proof.}
The necessity of the condition is clear.
Conversely, assume that the condition is satisfied.
It suffices to show that $N\cap U$
is constancy-enforcing for all $U$ in an open cover of~$M$.
We may therefore assume that $N$ is contained
in the domain $P$ of a ${{\mathbb K}}$-analytic diffeomorphism
$\psi\colon P\to Q$ onto an open subset $Q\subseteq{{\mathbb K}}^m$.
It then suffices to show that $\psi(N)$
is a constancy-enforcing subset of~$Q$
(or, equivalently, of~${{\mathbb K}}^m$).
We may therefore assume that $M={{\mathbb K}}^m$.
Let $Z$ be a ${{\mathbb K}}$-analytic
manifold and $f\colon Z\to {{\mathbb K}}^m$ be a ${{\mathbb K}}$-analytic
mapping such that $f(Z)\subseteq N$.
Given $x\in Z$, there exists a ${{\mathbb K}}$-analytic diffeomorphism
$\phi\colon U\to V$ from an open $x$-neighbourhood $U\subseteq Z$ onto an open
subset $V\subseteq{{\mathbb K}}^n$, with $\phi(x)=0$.
It suffices to show that $f\circ\phi^{-1}$ is constant on some
$0$-neighbourhood.
We may therefore assume that
$Z=B_\varepsilon(0)$ is an open ball in ${{\mathbb K}}^n$
for some norm $\|\cdot\|$ on~${{\mathbb K}}^n$ and some
$\varepsilon>0$. After a translation in the range, we may assume
that $f(0)=0$.
We may also assume that
\[
f(y)=\sum_{k=1}^\infty p_k(y)
\]
pointwise for $y\in B_\varepsilon(0)$
with certain homogeneous polynomials
$p_k\colon {{\mathbb K}}^n\to{{\mathbb K}}^m$
of degree~$k$.
For each $y\in B_\varepsilon(0)$, we obtain a ${{\mathbb K}}$-analytic map
\[
h\colon {{\mathbb D}}\to {{\mathbb K}}^m,\quad z\mapsto f(zy)=\sum_{k=1}^\infty z^k\,p_k(y)
\]
with $h({{\mathbb D}})\subseteq f(U)\subseteq N$.
By hypothesis, $h$ is constant on a $0$-neighbourhood
$V\subseteq {{\mathbb D}}$, and thus $h|_V=0$.
As a consequence, $p_k(y)=0$ for all $k\in{{\mathbb N}}$
(see second lemma in \cite[Part~II, Chapter~II]{Ser},
applied to power series in a single variable).
Thus $p_k=0$, as $y\in B_\varepsilon(0)$ was arbitrary
and $p_k$ is homogeneous. Hence $f=0$.
{\nopagebreak\hspace*{\fill}$\Box$\medskip\medskip\par}
\section{Proof of Theorem~1.1}
We shall deduce the theorem from the following
lemma. Let ${{\mathbb K}}:={{{\mathbb F}}\mkern+.2mu(\mkern-2.3mu(X)\mkern-2.3mu)}$ and
${{\mathbb O}}:={{\textstyle{{\mathbb F}}\mkern-1.7mu \left[\mkern-1mu \left[X\mkern+1mu\right]\mkern-1mu \right]}}$. Then ${{\mathbb D}}
:=\{z\in{{\mathbb K}}\colon|z|<1\}=X{{\mathbb O}}$.\\[2.5mm]
{\bf Lemma~3.1}
\emph{$H$ as in {\rm(\ref{theone})}
is a constancy-enforcing subset of ${{\mathbb K}}$.}\\[2.5mm]
{\bf Proof.}
By Lemma~2.6,
we have to show that each ${{\mathbb K}}$-analytic map
$f\colon {{\mathbb D}}\to{{\mathbb K}}$ with $f({{\mathbb D}})\subseteq H$
is constant on some $0$-neighbourhood.
After replacing $f$ with $z\mapsto f(\lambda z)$
for a small non-zero element $\lambda\in{{\mathbb D}}$,
we may assume that $f$
is given globally by a convergent power series,
\[
f(z)=\sum_{k=0}^\infty a_kz^k
\]
with $a_k\in {{\mathbb K}}$.
Then $a_0=f(0)\in H$. Since $H$ is a subgroup, it follows that
$f(y)-a_0\in H$ for all $y\in {{\mathbb D}}$.
We may therefore assume without loss of generality that $a_0=0$.
Pick $0\not=w\in{{\mathbb D}}$. Since the series for $f(w)$ converges,
the summands form a bounded sequence; thus
\[
M:=\sup\{|a_k|\,|w|^k\colon k\in{{\mathbb N}} \}<\infty.
\]
Choose $\theta\in {{\mathbb D}}$ such that $0<|\theta|\leq \frac{1}{M+1}$.
Then $|\theta^k|\leq \frac{1}{M+1}$ for all $k\in{{\mathbb N}}$ and thus
$|a_k(w\theta)^k|\leq 1$ for all $k\in{{\mathbb N}}$. After replacing $f$ with
$z\mapsto f(zw\theta)$, we may therefore assume that
the coefficients of $f$ satisfy
\[
|a_k|\leq 1\quad\mbox{for all $\,k\in{{\mathbb N}}$.}
\]
We show that $f=0$.
If $f\not=0$, then there exists a minimal $\ell\in {{\mathbb N}}$ such that $a_\ell\not=0$;
we shall derive a contradiction.
We can write
\[
a_\ell=bX^m
\]
with $m\in{{\mathbb N}}_0$ and $b\in {{\mathbb K}}$ such that $|b|=1$.
For all $n\in {{\mathbb N}}$ such that $n\geq m+1$,
\[
f(X^n)=a_\ell X^{\ell n}+R_n=b X^{\ell n+m}+R_n
\]
holds with
\[
R_n=\sum_{k=\ell+1}^\infty a_kX^{kn}\in X^{\ell n+m+1}{{\mathbb O}},
\]
since $a_k\in{{\mathbb O}}$ and $kn\geq(\ell+1)n=\ell n+ n\geq \ell n+m+1$.
Thus $f(X^n)\in H$ necessitates 
\[
\ell n+m\in\{2^j\colon j\in{{\mathbb N}}_0\},
\]
for all $n\geq m+1$. This is absurd.
{\nopagebreak\hspace*{\fill}$\Box$\medskip\medskip\par}

\noindent
{\bf Proof of Theorem~1.1.}
(a) Let $Z$ be a ${{\mathbb K}}$-analytic manifold and
$f\colon Z\to H$ be a map.
If $\iota_H\circ f$ is ${{\mathbb K}}$-analytic, then $f$ is locally constant,
since~$H$ is constancy-enforcing
(by Lemma~3.1).
Hence $f$ is ${{\mathbb K}}$-analytic as a map
to $H$, endowed with the discrete topology
(and corresponding ${{\mathbb K}}$-analytic manifold
structure modeled on ${{\mathbb K}}^0$).
If, conversely,
$f$ is ${{\mathbb K}}$-analytic as a map to $H$ with the discrete
topology, then~$f$ is locally constant,
whence also $\iota_H\circ f$ is locally constant and hence
${{\mathbb K}}$-analytic.\\[1mm]
(b) Endow $H$ with a ${{\mathbb K}}$-analytic manifold
structure which turns $\iota_H$ into a ${{\mathbb K}}$-analytic
map. Since $H$ is constancy-enforcing,
$\iota_H$ is locally constant.
Since $\iota_H$ is injective, we deduce that $H$
is discrete.
{\nopagebreak\hspace*{\fill}$\Box$\medskip\medskip\par}
\section{Proof of Theorem~1.2}
The point is that
${{\mathbb Z}}_p$ cannot be endowed with a ${{\mathbb K}}$-analytic
Lie group structure compatible with its natural topology,
for ${{\mathbb K}}={{{\mathbb F}}\mkern+.2mu(\mkern-2.3mu(X)\mkern-2.3mu)}$.
This can be deduced from more elaborate
general arguments (see \cite[Theorem~13.23]{Dix}
or further developments in~\cite{JZK}).
For the reader's convenience, we provide a direct proof.
Suppose that~${{\mathbb Z}}_p$ admits such a manifold structure,
modeled on ${{\mathbb K}}^n$ with $n\in{{\mathbb N}}$;
we shall derive a contradiction.
Write~$H$ for~${\mathbb Z}_p$ with this structure.
Let $q$ be the order of the finite
field~${{\mathbb F}}$; thus $q=p^m$ for some
$m\in{{\mathbb N}}$.
Choose the absolute value $|\cdot|$ on ${{\mathbb K}}$
such that $|X|=1/q$.
Let $\|\cdot\|$ be the maximum norm
on ${{\mathbb K}}^n$ given by $\|(z_1,\ldots, z_n)\|:=\max\{|z_1|,\ldots,|z_n|\}$
and define
\[
B_\varepsilon:=\{z\in {{\mathbb K}}^n\colon \|z\| <\varepsilon\},\quad \overline{B}_\varepsilon:=
\{z\in{{\mathbb K}}^n\colon \|z\|\leq \varepsilon\}
\]
for $\varepsilon>0$.
As ${{\mathbb D}}=X{{\mathbb O}}$, we have $B_1=X{{\mathbb O}}^n=X\mkern+1mu\overline{B}_1$.
More generally,
\[
B_{q^{-k}}=\overline{B}_{q^{-k-1}}=X^{k+1}\overline{B}_1\quad\mbox{for each integer~$k$.}
\]
Considering balls as compact, open subgroups
of $({{\mathbb K}}^n,+)$, we obtain the index
\[
\left[\mkern+2mu\overline{B}_1:X\overline{B}_1\right]=\left[{{\mathbb O}}:X{{\mathbb O}}\right]^n=q^n.
\]
Since ${{\mathbb K}}^n\to{{\mathbb K}}^n$, $y\mapsto X^ky$
is an automorphism of $({{\mathbb K}}^n,+)$, we deduce
that
$\left[\mkern+2mu\overline{B}_{q^{-k}}:\overline{B}_{q^{-k-1}}\right]=\left[X^k\overline{B}_1:
X^k(X\overline{B}_1)\right]
=q^n$ for each integer~$k$ and thus
\[
\left[\mkern+2mu\overline{B}_{q^{-k}}:\overline{B}_{q^{-\ell}}\right]=(q^n)^{\ell-k}\;\, \mbox{for all integers $\ell\geq k$.}
\]
By \cite[Part II, Chapter~IV, \S8, Theorem]{Ser},
$H$ has an open subgroup~$U$ which is a standard
group. Notably, there exists $r\in \; ]0,1]$ and an isomorphism
$\phi\colon U\to B_r$
of ${{\mathbb K}}$-analytic Lie groups
for a group multiplication $B_r\times B_r\to B_r$, $(x,y)\mapsto x*y$
on~$B_r$
with neutral element~$0$, such that
$B_\varepsilon$ is a normal subgroup of $B_r$ for all $\varepsilon\in\,]0,r]$
(cf.\ \S9 in the book chapter just cited).
Moreover, the cosets of $B_\varepsilon$
in $(B_r,*)$ and $(B_r,+)$
coincide, for all $\varepsilon\in \,]0,r]$
(see Theorem~1 in \cite[Part~II, Chapter~IV, \S9]{Ser}).
After shrinking~$r$ if necessary, we may assume
that the components $f_1,\ldots, f_n\colon B_r\to{{\mathbb K}}$
of the $p$th
power map
\[
f\colon B_r\to B_r,\quad z\mapsto z^p=z*\cdots * z
\]
in $(B_r,*)$
are given by power series
\[
f_j(z)=\sum_{\alpha\in({{\mathbb N}}_0)^n}a_{\alpha,j} z^\alpha
\]
on $B_r$
for all $j\in\{1,\ldots, n\}$,
with $\sum_{\alpha\in ({{\mathbb N}}_0)^n}|a_{\alpha,j}|\,r^{|\alpha|}<\infty$.
We use multi-index notation here, as in~\cite{Ser};
notably, $|\alpha|:=\alpha_1+\cdots+\alpha_n$
for $\alpha=(\alpha_1,\ldots,\alpha_n)\in({{\mathbb N}}_0)^n$,
$z^\alpha:=z_1^{\alpha_1}\cdots z_n^{\alpha_n}\in{{\mathbb K}}$
for $z=(z_1,\ldots, z_n)\in{{\mathbb K}}^n$,
and $a_{\alpha,j}\in{{\mathbb K}}$.
Setting $a_\alpha:=(a_{\alpha,1},\ldots,a_{\alpha,n})\in{{\mathbb K}}^n$,
we have
$M:=\sum_{\alpha\in({{\mathbb N}}_0)^n}\|a_\alpha\|\, r^{|\alpha|}<\infty$
and
\[
f(z)=\sum_{\alpha\in ({{\mathbb N}}_0)^n}a_\alpha z^\alpha\quad\mbox{in $\,{{\mathbb K}}^n\,$ for all $\, z\in B_r$,}
\]
as the limit of the net of finite partial sums.
Then $a_0=0$ as $f(0)=0^p=0$ and $a_\alpha=0$ for all
multi-indices $\alpha$ of length $|\alpha|=1$.
In fact, the power map~$f$ has derivative
$Df(0)=p\hspace*{.4mm}\mbox{id}_{{{\mathbb K}}^n}
=0\in \mbox{End}_{{\mathbb K}}
({{\mathbb K}}^n)$ at~$0$
since ${{\mathbb K}}$ has characteristic~$p$
(see \cite[Part~II, Chapter~IV, \S7]{Ser}).
For all $z\in B_r$ and $\alpha\in ({{\mathbb N}}_0)^n$ with $|\alpha|\geq 2$,
\[
|z^\alpha|\leq (\|z\|)^{|\alpha|}=\left(\frac{\|z\|}{r}\right)^{|\alpha|}
r^{|\alpha|}\leq \left(\frac{\|z\|}{r}\right)^2r^{|\alpha|}
\]
holds, whence
\[
\|f(z)\|\leq \sum_{|\alpha|\geq 2}\|a_\alpha\|\, |z^\alpha|
\leq \left(\frac{\|z\|}{r}\right)^2\sum_{|\alpha|\geq 2}\|a_\alpha\| r^{|\alpha|}
=C \|z\|^2
\]
with $C:=M/r^2$.
Thus
\[
f(B_\varepsilon)\subseteq B_{C\varepsilon^2}
\]
for all $\varepsilon\in\;]0,r]$.
Notably, $f(B_\varepsilon)\subseteq B_\varepsilon$
if $\varepsilon\in \,]0,r]$ such that $C\varepsilon\leq 1$.
Choose $\ell\in{{\mathbb N}}$ so large that
$Cq^{-\ell}\leq q^{-2}$ and $q^{-\ell}\leq r$.
Then $C(q^{-\ell})^2=Cq^{-\ell}q^{-\ell}\leq q^{-\ell-2}$
and thus, choosing $\varepsilon:=q^{-\ell}$,
\[
f\left(\overline{B}_{q^{-\ell-1}}\right)=f\left(B_{q^{-\ell}}\right)\subseteq B_{C(q^{-\ell})^2}\subseteq B_{q^{-\ell-2}}
=\overline{B}_{q^{-\ell-3}},
\]
whence
\[
\left[\mkern+2mu\overline{B}_{q^{-\ell-1}}:f(\overline{B}_{q^{-\ell-1}})\right]
\geq \left[\mkern+2mu \overline{B}_{q^{-\ell-1}}:\overline{B}_{q^{-\ell-3}}\right]
=(q^n)^2.
\]
But $W:=\phi^{-1}(\overline{B}_{q^{-\ell-1}})$
is a closed subgroup of~${{\mathbb Z}}_p$ and hence
of the form $p^k{{\mathbb Z}}_p$ for some $k\in{{\mathbb N}}_0$.
Since $\phi$ is a homomorphism of groups and~$f$ the $p$th power map,
we have $\phi(pW)=f(\phi(W))$.
Thus
\[
p=[p^k{{\mathbb Z}}_p:p^{k+1}{{\mathbb Z}}_p]=[W,pW]=[\phi(W),f(\phi(W))]\geq (q^n)^2\geq p^2,
\]
contradiction.
{\nopagebreak\hspace*{\fill}$\Box$\medskip\medskip\par}
\section{Proof of Theorem~1.4}
Give ${{\mathbb K}}:={{{\mathbb F}}\mkern+.2mu(\mkern-2.3mu(X)\mkern-2.3mu)}$ the absolute
value taking $X$ to $1/p$.
By Lemma~2.6,
we have to show that each ${{\mathbb K}}$-analytic map
$f\colon {{\mathbb D}}\to{{\mathbb K}}$ with $f({{\mathbb D}})\subseteq H:={{\textstyle{{\mathbb F}}\mkern-2.2mu \left[\mkern-2.8mu \left[X^\ell\mkern+1.5mu\right]\mkern-2.8mu \right]}}$
is constant on some $0$-neighbourhood.
As in the proof of Lemma~3.1,
we may assume that
$f(0)=0$ and that,
for suitable $a_k\in {{\mathbb K}}$ with $|a_k|\leq 1$,
\[
(\forall z\in{{\mathbb D}})\qquad f(z)=\sum_{k=1}^\infty a_k\, z^k.
\]
We show that $f=0$.
If $f\not=0$,
we get a contradiction:
First, note that $f'(z)=0$
for all $z\in{{\mathbb D}}$.
In fact, if $f'(z_0)\not=0$,
the inverse function theorem
would show that $f({{\mathbb D}})$
is a neighbourhood of $f(z_0)$ in ${{\mathbb K}}$;
this is absurd, as $f({{\mathbb D}})$ is a subset of $H$
and $H={{\textstyle{{\mathbb F}}\mkern-2.2mu \left[\mkern-2.8mu \left[X^\ell\mkern+1.5mu\right]\mkern-2.8mu \right]}}$
has empty interior in~${{\mathbb K}}$.
Hence
\[
(\forall z\in{{\mathbb D}})\qquad
\sum_{k=1}^\infty k\,a_k\, z^{k-1}=0,
\]
whence $a_k=0$ for all $k\in{{\mathbb N}}$ such that $k\not\in p{{\mathbb N}}$
(see the second lemma in \cite[Part~II, Chapter~II]{Ser}).
As we suppose that $f\not=0$,
some coefficient of~$f$ is non-zero.
Hence, there exists a maximal
$n\in{{\mathbb N}}$ such that $a_k=0$
for all $k\in{{\mathbb N}}$ such that $k\not\in p^n {{\mathbb N}}$.
Abbreviate $q:=p^n$.
By maximality of $n$, there exists $m\in{{\mathbb N}}$
with $a_m\not=0$
such that $m\in q{{\mathbb N}}$ but $m\not\in pq{{\mathbb N}}$
and thus $m=q\mu$ with $\mu\in {{\mathbb N}}\setminus p{{\mathbb N}}$.
By the preceding,
\[
(\forall z\in{{\mathbb D}})\;\; f(z)=\sum_{k=1}^\infty a_{qk}\,z^{qk}.
\]
We now consider ${{\mathbb K}}$
as a subfield of ${{\mathbb L}}:={{{\mathbb F}}\mkern+.2mu(\mkern-2.3mu(\mkern+.4mu{}Y\mkern+.2mu )\mkern-2.3mu)} \cong {{\mathbb K}}$
using the embedding
\[
\sum_{j=j_0}^\infty c_jX^j\mapsto \sum_{j=j_0}^\infty c_j Y^{qj}
\]
for $j_0\in{{\mathbb Z}}$, $c_j\in{{\mathbb F}}$.
Thus ${{\mathbb K}}=({{\mathbb L}})^q$ is the image
of ${{\mathbb L}}$ under the $n$th power of the
Frobenius endomorphism of~${{\mathbb L}}$.
We give ${{\mathbb L}}$ the absolute value
(also denoted $|\cdot|$)
which extends the given absolute value
$|\cdot|$ on~${{\mathbb K}}$;
thus $|Y|=\sqrt[q]{1/p}$.
For each $k\in{{\mathbb N}}$, there exists $b_k\in{{\mathbb L}}$
such that $(b_k)^q=a_{qk}$.
Then $b_\mu\not=0$ and
\[
g\colon W \to {{\mathbb L}},\quad w\mapsto\sum_{k=1}^\infty b_k\mkern+1mu w^k
\]
is an ${{\mathbb L}}$-analytic function
on $W:=\{w\in{{\mathbb L}}\colon |w|<1\}$
as all $|b_k|\leq 1$.
Now
\[
f(z)=(g(z))^q\quad\mbox{for all $z\in {{\mathbb D}}$,}
\]
by construction.
Since $\mu \mkern+1mu b_\mu\not=0$ holds in~${{\mathbb L}}$,
the function $g'\colon W\to {{\mathbb L}}$,\linebreak
$w\mapsto \sum_{k=1}^\infty k\mkern+1mu b_k\, w^{k-1}$
does not vanish identically.
Let $\nu\in{{\mathbb N}}_0$ be minimal such that $(\nu+1) b_{\nu+1}\not=0$;
then\vspace{-.4mm}
\[
g'(w)=w^\nu h(w)
\]
with $h(w):=\sum_{k=\nu+1}^\infty k\mkern+1mu b_k\mkern+1mu w^{k-1-\nu}$
and $h(0)\not=0$.
There exists $r\in\,]0,1]$
such that $h(w)\not=0$
for all $w\in {{\mathbb L}}$ such that $|w|<r$.
There exists $z_0\in {{\mathbb K}}\setminus \{0\}$
such that $|z_0|<r$; then $g'(z_0)=(z_0)^\nu h(z_0)\not=0$.
Now
\[
g(w)=g(z_0)+g'(z_0)(w-z_0)+(w-z_0)R(w)
\]
for a continuous function $R\colon W\to{{\mathbb L}}$ such that
$R(z_0)=0$. Since $g'(z_0)=\sum_{k=1}^\infty k\mkern+1mu b_k\,(z_0)^{k-1}$
with $|k\mkern+1mu b_k|\leq 1$, we have $|g'(z_0)|\leq 1$.
Thus, we can write
$g'(z_0)=a Y^\tau$ with $\tau\in {{\mathbb N}}_0$
and $a\in{{\mathbb L}}$ such that $|a|=1$.
Then $a=b+c$ with $b\in{{\mathbb F}}^\times$ and $c\in Y\mkern+1mu{{\textstyle{{\mathbb F}}\mkern-1.7mu \left[\mkern-1mu \left[\mkern+1mu Y\mkern+1mu\right]\mkern-1mu \right]}}$.
There is $s\in\,]0,1]$ such that
$|R(w)|<|g'(z_0)|$ and
thus\vspace{-.4mm}
\[
g'(z_0)+R(w) =\mkern+1mu b \mkern+1mu Y^\tau+\,\mbox{higher order terms}
\]
for all $w\in W$ such that $|w-z_0|<s$.
We choose $j\in{{\mathbb N}}$ so large that $(1/p)^j<s$.
Then $|X^j|<s$ and $|X^{j+1}|<s$.
Moreover,
\[
H\ni f(z_0+X^j)-f(z_0)=
(g(z_0+X^j)-g(z_0))^q= b^q X^{\tau+qj}+\,\mbox{higher order terms},
\]
whence $\tau+ qj\in \ell \mkern+1mu {{\mathbb N}}$. Likewise,\vspace{.3mm}
$f(z_0+X^{j+1})-f(z_0) \in H$ implies that $\tau+q(j+1)\in\ell\mkern+1mu {{\mathbb N}}$.
Hence $q\in \ell\mkern+1mu {{\mathbb N}}$, which is absurd. $\,\square$\vspace{-4mm}
{\footnotesize
{\bf Helge  Gl\"{o}ckner}, Universit\"at Paderborn, Institut f\"{u}r Mathematik,\\
Warburger Str.\ 100, 33098 Paderborn, Germany;
{\tt glockner\symbol{'100}math.upb.de}}
\end{document}